\documentclass[11pt]{article}
\usepackage{hyperref}
\usepackage{color}
\usepackage[square,authoryear]{natbib}
\usepackage{marsden_article}

\begin{document}

\AllowFootnote
\SetVersionDate{May 2, 2000}
\editedby{jem}


\title{The Anisotropic Averaged Euler Equations}
\author{Jerrold E. Marsden \\
Control and Dynamical Systems 107-81 \\
California Institute
of Technology \\
Pasadena, CA 91125 \\
{\footnotesize marsden@cds.caltech.edu}
\and
Steve Shkoller \\
Department of Mathematics\\
University of California \\
Davis, CA 95616\\
{\footnotesize shkoller@math.ucdavis.edu}}
\date{September, 1999; this version: May 3, 2000\\[10pt]
{\it Dedicated to Stuart Antman on the occasion of his 60th birthday}}


\maketitle

\tableofcontents

\begin{abstract}
The purpose of this paper is to derive the anisotropic averaged Euler
equations and to study their geometric and analytic properties.  These new
equations involve the evolution of a mean velocity field and an advected
symmetric tensor that captures the fluctuation effects. Besides the
derivation of these equations, the new results in the paper are smoothness
properties of the equations in material representation, which gives
well-posedness of the equations, and the derivation of a corrector to the
macroscopic velocity field. The numerical implementation and physical
implications of this set of equations will be explored in other
publications.
\end{abstract}

\section{Introduction}
A fundamental problem in turbulent fluid dynamics is the difficulty in
resolving the many spatial scales that are activated by the
complicated nonlinear interactions. It is a challenge to produce models
that capture the large scale flow, while correctly modeling the
influence of the small scale dynamics. While there are many efforts in this
direction, the goal of the present paper is to introduce a new method that
is based on the combination of two basic ideas: the use of an
ensemble averaging
that represents a spatial sampling of material particles over small
spatial scales, and the use  of asymptotic expansions together with this
averaging on the level of the variational principle.

Our approach is conceptually similar to the method of Reynolds averaging and
Large Eddy Simulation techniques, but has the advantage of 1) not needing
an additional closure model and 2) automatically providing a small scale
corrector to the macroscopic flow field.

Our methodology has some interesting connections with the method of Optimal
Prediction introduced by \cite{ChKaKu1999}, which will be explored in
future publications.

In the body of the paper we shall comment on a comparison between our
approach and that of \cite{ChFoHoOlTiWy1998} and \cite{Holm1999}, which
produces different equations.

\subsection{A Brief Review of the Euler and Isotropic Averaged Euler
Equations} \label{review-section}

\paragraph{A Brief History.} There has been much recent interest in the
averaged Euler equations for ideal fluid flow. In this paper we will focus
on the geometry and analysis of a related set of equations, which we call
the {\it anisotropic averaged Euler equations.} The original averaged Euler
equations appear as a special {\it isotropic} case of the more general
equations.

The isotropic averaged Euler equations on all of ${\mathbb R}^n$ first
appeared in the context of an approximation to the Euler equations in
\cite{HoMaRa1998a} and some of its variational structure was developed in
\cite{HoMaRa1998b};  this variational structure retains the quadratic form
of the variational structure for the original Euler equations, so that the
equations can be viewed as describing a certain geodesic flow in a sense
similar to the work of \cite{Arnold1966} and \cite{EbMa1970}.

Remarkably, these equations are mathematically identical to the well-known
inviscid  {\it second grade fluids} equations introduced by
\cite{RiEr1955}. The geometric analysis of these equations, including
well-posedness and other  analytic properties, was developed in
\citeauthor{Shkoller1998} [\citeyr{Shkoller1998},
\citeyr{Shkoller2000}] and \cite{MaRaSh2000}.  These references also
discuss the relation to the second-grade fluid literature.

In \cite{OlSh2000}, the link with the vortex blob method was established;
therein, it was shown that the vortex blob numerical algorithm generates
unique global weak solutions to the averaged Euler equations.   These weak
solutions induce a {\it weak} coadjoint action on the vector space of
vorticity functions, modeled as the space of Radon measures.  The existence
of such a weak coadjoint action makes rigorous the formal constructions of
\cite{MaWe1983} on the geometry of point-vortex and vortex blob dynamics.

The works of \cite{ChFoHoOlTiWy1998} and \cite{Holm1999} formulated
equations for the slow time dynamics of fluid motion by averaging over fast
time fluctuations about the mean;  that approach, founded on a Reynolds
decomposition translated over the Lagrangian parcel, and the resulting
system of equations, is different from  our approach and from the results
that we shall present. We give a few more details on the comparison in the
body of the paper.

\paragraph{The Euler Equations as Geodesics and Notation.}

It is well-known how to view the Euler equations as geodesics on the
group of diffeomorphisms and that this view has concrete analytical
advantages, due to the work of \cite{Arnold1966} and \cite{EbMa1970}. In
particular, this work shows that the equations define a smooth vector field
(a spray) on the group of diffeomorphisms, that is, in Lagrangian (or
material) representation. The reduction of the equations from material to
spatial (Eulerian) representation may be viewed by the classical and
general technique of Euler-Poincar\'e reduction (see \cite{MaRa1999} and
\cite{HoMaRa1998b} for an exposition and further references) and this view
is a helpful guide to understanding other fluid theories as well.

The geometric view of fluid mechanics, along with a careful understanding of
the averaging process, will be basic to the present paper, so we briefly
review the salient features of the theory  for the reader's convenience,
and to establish notation.

Let $(M, g )$ be a $ C ^{\infty}$ compact, oriented $n$-dimensional
Riemannian manifold with $C ^{\infty}$ boundary (possibly empty). Of
course open regions with smooth boundary in the plane or space are key
examples. The Riemannian volume form associated with the metric $g$ is
denoted $\mu$.

The Euler equations for the velocity field $u$ of an ideal,
incompressible, homogeneous fluid moving on
$M$ (such as a region $\Omega$ in $\mathbb{R}^2$ or
$\mathbb{R}^3$) are
\begin{equation}
\label{eef} \frac{\partial u}{\partial t}
+ (u \cdot \nabla )u  = -\nabla p
\end{equation}
with the constraint $\operatorname{div}\, u = 0$ and the
boundary condition that $u$ is tangent to the boundary,
$\partial M.$ The pressure $p$ is determined by the incompressibility
constraint. The nonlinear term $(u \cdot \nabla )u$ is
interpreted in the context of manifolds to be $\nabla _{u}
u$, the covariant derivative of $u$ along $u$.
In Euclidean coordinates, these equations are given as follows (using the
summation convention for repeated indices):
\[
\frac{\partial u^i}{\partial t}
+ u^j
\frac{\partial u^i }{\partial x^j}
= - \frac{\partial p}{\partial x^i}\, ,
\]
and on a Riemannian manifold (or in curvilinear coordinates in Euclidean
space), the Euler equations take the following coordinate form:
\[
\frac{\partial u^i}{\partial t}
+ u^j \frac{\partial u^i }{\partial x^j}
+ \Gamma^i_{jk} u ^j u ^k
= - g^{ij}\frac{\partial p}{\partial x ^j}\,,
\]
where $g _{ij}$ are the components of the Riemannian metric $g$,
$g ^{ij} = [g_{ij}]^{-1}$,
and $\Gamma^i_{jk}$ are the associated Christoffel symbols. Using
covariant derivative notation, these coordinate equations read
\[
\frac{\partial u^i}{\partial t}
+ u^j  u^i _{\; ; j}
= - g^{ij}p _{,j}.
\]
We let the flow of the time dependent vector field $u(t,x)$ be denoted by
$\eta(t,x)$ so that
$$
\frac{\partial}{\partial t} \eta(t,x)  = u (t,\eta(t,x)),
$$
with $\eta(0,x) = x$ for all $x$ in $M$. For each $t$, we denote the map
$\eta(t, \cdot)$ by $\eta_t$ so that $\eta_0=e$,  the identity map.  Thus,
the map $x \mapsto \eta _t (x)$ gives the particle placement field for the
fluid. Because of the incompressibility, each map
$\eta_t$ is volume preserving and is a diffeomorphism.

We shall be working with vector fields $u$ of Sobolev class $H ^s$ for $s >
(n/2) + 1$ and, correspondingly, $\eta _t \in \mathcal{D}_{\mu}^s$, the
group of $H ^s$-volume preserving diffeomorphisms. If there is any danger
of confusion, we shall write $\mathcal{D}_{\mu}^s(M)$ to indicate the
underlying manifold $M$.  See \cite{EbMa1970} and \cite{Shkoller2000} for
some basic properties of Hilbert class diffeomorphism groups for manifolds
with boundary.

Arnold's theorem on the Euler equations may be stated as follows: {\it A
time dependent velocity field $u $ satisfies the Euler equations iff the
curve $\eta _t $ is a geodesic of the right invariant $L^2$-metric on
$\mathcal{D}_{\mu}^s$.}

This $L^2$-metric is defined as follows. The tangent space to
$\mathcal{D}_{\mu}^s$ at the identity is identified with the space
$\mathfrak{X}_{\operatorname{div}}^s$, the space of $H^s$ divergence free
vector fields on $M$ that are tangent to the boundary $\partial M $. The
right invariant
$L^2$-metric is defined to be the weak Riemannian metric on
$\mathcal{D}_{\mu}^s$ whose value at the identity is
\[
\left\langle u, w \right\rangle_{L^2} = \int _M \left\langle u (x), w (x)
\right\rangle_x \mu(x),
\]
where we write the pointwise inner product as
$ \left\langle u (x), w (x) \right\rangle_x =
g (x) (u (x), w (x) )$, and the pointwise norm $|u(x)|^2 = \left\langle
u(x),u(x)\right\rangle_x$.

As we shall explain shortly, with the maturation of Euler-Poincar\'e theory,
Arnold's theorem becomes an easy consequence of more general and rather
simple results.

\paragraph{Lie Derivative and Vorticity Form.} As is well-known, the Euler
equations can be written  in terms of Lie derivatives as
\begin{equation} \label{Lie_Euler.equation}
\frac{\partial u ^{\flat} }{\partial t}  + \pounds _u u
^{\flat} =  \mathbf{d} \left( \frac{1}{2} | u | ^2 - p \right) =
- \mathbf{d} p ^\prime,
\end{equation}
where $u ^\flat$ is the one-form associated to the vector field $u$ via
the metric, and $\pounds _u u ^\flat$ denotes the Lie derivative of the
one-form $u ^\flat$ along $u$. Taking the exterior derivative of
(\ref{Lie_Euler.equation}) gives the familiar advection equation for
vorticity:
$$ \frac{\partial \omega}{\partial t} + \pounds _{u } \omega = 0, $$
where $\omega = \mathbf{d} u ^\flat$ is the vorticity, thought of as a
two-form.   In 2D, $\omega$ is identified with a scalar and is traditionally
thought of as the 2D-curl of the velocity field, while in 3D, $\omega$
may be identified (using the volume-form $\mu$) with a vector field which
is traditionally obtained by taking the curl of $u$.

The vorticity equation is the infinitesimal version of the following
advection property:
$$\omega _t = (\eta _t) _\ast \omega _0.$$
Of course in
two dimensions, this gives the usual advection of vorticity as a
function, while in three  (or higher) dimensions, the advection is
understood in terms of advection of two-forms.

The Euler equations have both an interesting Hamiltonian structure
in terms of Poisson brackets (a Lie-Poisson bracket) and a variational
structure. In this paper we shall be working primarily with the
variational structure; the Hamiltonian structure, along with references
to the literature may be found in \cite{MaWe1983}, \cite{ArKh1998} and
\cite{MaRa1999}.

\paragraph{Lagrangian and Variational Form.} The Lagrangian is given by the
total kinetic energy of the fluid; in spatial representation, this
Lagrangian is
\begin{equation}  \label{lagrangian_euler.equation}
L(u) = \frac{1}{2} \int _M | u(x) | ^2 \mu.
\end{equation}
The corresponding (unreduced) La\-gran\-gian on $T \mathcal{D}_{\mu}^s$ is
given by
\begin{equation} \label{Lagrangian_euler.equation}
\mathcal{L} ( \eta, \dot{\eta} ) = \frac{1}{2} \int _M
g(\eta(x)) (\dot\eta(x),\dot\eta(x)) \mu.
\end{equation}
Hamilton's principle on $\mathcal{D}_{\mu}^s$ applied to the Lagrangian
$\mathcal{L}$ gives geodesics on this group. Euler-Poincar\'e reduction
techniques (see \cite{MaRa1999}) show that this variational
principle reduces to the following principle in terms of Eulerian
velocities:
$$
\delta \int _a ^b L (u) \,dt = 0,
$$
which should hold for all variations
$\delta u$ of the form
$$
        \delta u = \dot{w} + [u, w],
$$
where $w$ is a time dependent vector field (representing the infinitesimal
particle displacement) vanishing at the temporal endpoints\footnote{The
constraints on the allowed variations of the fluid velocity field are
commonly known as ``Lin constraints''. This itself has an interesting
history, going back to Ehrenfest, Boltzmann,  Clebsch, Newcomb and
Bretherton, but there was little if any contact with the heritage of Lie
and Poincar\'{e} on the subject.}. Here, $[w, u]$ denotes the usual
Jacobi--Lie bracket of vector fields. One readily checks that this reduced
principle yields the standard Euler equations. This simple computation is
the heart of Arnold's theorem.

\paragraph{Analytical Issues.} While the Eulerian (spatial) representation
has been emphasized in most analytic studies of the Euler equations,  fluid
motion viewed on the Lagrangian (material) side has some distinct
advantages.  For example, \cite{EbMa1970} proved that the flow, solving the
Euler equations, on the volume-preserving diffeomorphism  group
$\mathcal{D}_\mu^s$, $s>n/2+1$, is smooth in time.  They derived a number of
interesting consequences from this result, including theorems on the
convergence of solutions of the Navier-Stokes equations to solutions of the
Euler equations as the viscosity goes to zero when $M$ has no boundary. In
addition, \cite{MaPu1994} analyzed the Lagrangian flow map to establish
sharp well-posedness of the 2D Euler equations and prove convergence of the
vortex blob algorithm.  In many cases, the Lagrangian framework is, in
fact, the more natural setting to study the behavior of solutions, and we
shall emphasize this point of view.

\subsection{The Averaged Euler Equations}
\paragraph{The Isotropic Averaged Euler Equations.} Let $\alpha$ be a
positive constant. In Euclidean space and in Euclidean coordinates, the
isotropic averaged Euler equations (inviscid second-grade fluids
equations)\footnote{These are also known as the Euler-$\alpha$ equations.}
read:
$$
\frac{\partial v^i}{\partial t}
+ u^j \frac{\partial v^i}{\partial x^j}
-\alpha^2
\left[
\frac{\partial u^j}{\partial x^i}
\right]
\Delta  u_j = - \frac{\partial p}{\partial x ^i}\,,
$$
where $v = u - \alpha^2\Delta  u$ and $\Delta$ denotes the
componentwise Laplacian, and there is a summation over repeated indices
(in Euclidean coordinates, as is common, we make no distinction between
indices up or down). While there are several choices, the no slip boundary
conditions $u = 0$ are often used for this model.

\paragraph{Rate of Deformation Tensor.} One of the interesting things that
comes out of a careful derivation of the equations is the natural occurence
of the {\bfi rate of deformation tensor}, which is defined by
$$
\operatorname{Def} u = \frac{1}{2} \left(\nabla u + (\nabla u ) ^T \right)
$$
which we write in coordinates as:
$$
  \left( \operatorname{Def} u \right)^i_{\; j} =
\frac{1}{2} \left(u^i_{\; ;j} + u^j_{\; ;i} \right).
$$
We also let $\operatorname{Def}u^\flat =
\frac{1}{2} \left[\nabla u^\flat + (\nabla u^\flat ) ^T \right]$ which we
write in coordinates as
$$D _{ij} = \big( \operatorname{Def} u^\flat \big)_{ij} =
\frac{1}{2} \big(u_{i;j} + u_{j;i} \big).$$
Note that this is exactly the Lie derivative of the metric tensor; that is,
$\operatorname{Def} u^\flat = \pounds _{ u } g $, which is sometimes called
the {\it Killing tensor}.

\paragraph{Smoothness Properties.} Results on smoothness
of the La\-gran\-gian flow map for the averaged Euler equations were given
in \cite{Shkoller1998} on compact boundaryless Riemannian manifolds, and in
\cite{MaRaSh2000} on compact Euclidean domains.  The problem of how to
formulate this system on compact Riemannian manifolds {\it with boundary}
was solved in \cite{Shkoller2000};  the equations take the form
$$
\partial_t(1-\alpha^2\Delta _r)u + \nabla_u(1-\alpha^2\Delta _r)u
-\alpha^2 (\nabla u)^t \cdot \Delta _r u   = -\operatorname{grad} p,
$$
together with the constraint $\operatorname{div} u = 0,$ and
with appropriate initial conditions $u(0)   =u_0$, as well as boundary
conditions.  The symbol $\Delta_r$ is the operator $\operatorname{Def}
^*\operatorname{Def} $ acting on divergence-free vector fields, where
$\operatorname{Def} ^*$ is the $L^2$ formal adjoint of the (rate of)
deformation operator $\operatorname{Def} $.  Explicitly,
\begin{equation} \label{laplacian.equation}
\Delta _r = -(d \delta + \delta d) + 2\operatorname{Ric}.
\end{equation}
As with the usual Euler equations, the function $p$ is
determined from the incompressibility condition.

\paragraph{Lie Derivative Form---The Isotropic Equations.} The averaged
Euler equation can be neatly written in terms of Lie derivatives:
\begin{equation}
\partial _t v ^\flat + \pounds _u v ^\flat = -  d p ,
\end{equation}
where $v ^\flat = ( 1 -\alpha ^2 \Delta _r ) u ^\flat$.
\paragraph{The Anisotropic Averaged Euler Equations.} These equations,
which are the main subject of the present paper, and which will be derived
in \S\ref{derivation_section}, will now be stated. The basic variables that
are evolving in the anisotropic averaged Euler equations are the
{\bfi macroscopic velocity
field $u$} and a symmetric tensor field $F$ on $M$; the tensor field $F$
will be interpreted as the {\bfi contravariant spatial fluctuation tensor}
and it
will keep track of the anisotropy of the fluid deviations from the macroscopic
flow. These equations also depend on a choice of length scale $\alpha$.

  It is convenient to define the linear operator
${\mathcal C}:
\mathfrak{X}_{\operatorname{div}}^s
\cap H^1_0 \rightarrow H^{s-2}$, $s\ge 1$, by
$$\mathcal{C}u := \operatorname{Div}\left[ C :\nabla
u^\flat\right],$$
where $\flat$ is the map from vector fields  to one-forms associated with
the metric $g$, and the fourth-rank symmetric positive tensor $C$ is the 
symmetrization of the tensor $F\otimes g^{-1}$, given in local coordinates
by
$$ C^{ijkl}= {\frac{1}{4}}\left(
F^{lj}g^{ik}+ f^{kj}g^{il}+F^{li}g^{jk}+F^{ki}g^{jl}
\right).  $$

With this notation, the {\bfi anisotropic averaged Euler equations} on
 manifolds are
\begin{align}
\partial_t (1-\alpha^2\mathcal{C})u &+ \nabla_u(1-\alpha^2\mathcal{C})u
- \alpha^2
[\nabla u]^t \cdot \mathcal{C}u +
2\alpha^2 F : \left[\nabla (\operatorname{Def} u^\flat) ^2
\right]^\sharp \nonumber\\
& -
4\alpha^2 \operatorname{Div}\left( \left(\operatorname{Def}
u\right)^2 \cdot F \right)
= -\operatorname{grad} p,
\nonumber
\end{align}
together with the {\bfi advection equation}
$$ \partial_t F + \pounds_uF =0, $$
the incompressibility constraint
$\text{div }u=0$, initial data $u(0)=u_0$ and $F(0)=F_0$,
and the Dirichlet boundary condition $u=0$.

\paragraph{Lie Derivative Form---Anisotropic Equations.} The anisotropic
averaged Euler equations can also be written using Lie derivatives as
\begin{equation}
\partial _t v^\flat + \pounds _u v^\flat  + \left[
2\alpha^2 F : \left[\nabla (\operatorname{Def} u^\flat) ^2
\right]^\sharp
   -
4\alpha^2 \operatorname{Div}\left( \left(\operatorname{Def}
u\right)^2 \cdot F \right)
\right] ^\flat = -  dp,
\label{Lie_average.equation}
\end{equation}
where $v^\flat = (1 -\alpha ^2 \mathcal{C}^\flat) u^\flat$,
where ${\mathcal C}^\flat u^\flat = \text{Div}[C: \nabla u^\flat]^\flat$.

\paragraph{Coordinate Form.} In local coordinates,
the anisotropic averaged Euler equations become
\begin{align*}
& \partial_t \left(u^i - \alpha^2[C^{ijkl}u_k,_j],_l \right)
+\left( u^i - \alpha^2[C^{ijkl} u_k,_j],_l\right),_m u^m
-\alpha^2 u_m,_i[C^{mjkl}u_k,_j],_l \\
& \qquad  + 2\alpha^2 F^{kj} [D_{km}g^{mn}
D_{nj}],_i
      -4\alpha^2[F^{kj} D_{im}g^{mn} D_{nj}],_k   \\
& \qquad \quad = -p,_i
\end{align*}
together with the advection equation
\[
\partial_t F^{ij} + F^{ij},_k u^k - F^{kj}u^i,_k - F^{ik} u^j,_k = 0,
\]
with the constraint $u ^i_{\; , i } = 0 $, given initial conditions $u ^i
(0) = u ^i _0$,  and with the no-slip conditions $u ^i = 0$ on the
boundary.  If the metric $g$ is not the Euclidean metric $\delta_{ij}$, then
the partial derivatives above should be interpreted as arising from
the Levi-Civita covariant derivative associated to $g$.

\subsection{Outline of the Main Results.} The main results of the present work are as follows:
\begin{enumerate}
\item We derive, in a systematic way, the first order averaged Lagrangian
given in coordinates by
$$
L^\alpha_1(u,F) = \frac{1}{2} \int_M \left\{
g_{ik}u^i u^k + 2 \alpha^2 g^{ik} F^{jl} D _{ij}
D _{kl }  \right\}
[\text{det } g]^{\frac{1}{2}} dx.
$$
and, using the calculus of variations, derive the associated anisotropic
averaged Euler equations as the corresponding Euler-Poincar\'e equations.
The Euler-Poincar\'e technique was also used in \cite{Holm1999}, but the
Lagrangian and associated equations are different. In particular, the
principles and philosophy governing the derivation of the Lagrangian are
completely  different.
\item We show that the equations are well posed; in fact, we show more,
namely that the corresponding Lagrangian flow map is smooth in time in the
appropriate Sobolev topology.
\item Another important achievement is that while the macroscopic velocity
field $u$ is computed on spatial scales larger than $\alpha$, we are able
to obtain a corrector for this macroscopic field to order $\alpha^2$. This
is done in \S\ref{corrector_subsection} and is similar to what one does in
the theory of homogenization.
\end{enumerate}

\section{The Derivation}

\subsection{Introduction}
This section presents a new method for constructing
models of hydrodynamics which takes into account the fundamental
idea that a {\it fluid particle} is not a point, but rather a
{\it collection of points} forming a {\it representative sampling}.
Our approach is founded upon a certain type of Lagrangian ensemble averaging
performed at the level of the variational principle.
A similar idea on the level of the equation
itself, as opposed to the variational principle, was used by
\cite{BaPr1993} for deriving models of damage propagation.

\paragraph{Naive Averaging Does not Work.} We first explain why the naive
approach to spatially averaging a quadratic Lagrangian or Hamiltonian
does not suffice.  As a simple example, consider the Lagrangian on scalar
functions on $\mathbb{R}^n$ given by
$L(u)= \frac{1}{2} \int_{\mathbb{R}^n} u^2(x) dx $ and for a given positive
constant
$\alpha$, define a new averaged Lagrangian by
\[
L^\alpha(u) = \int_{\mathbb{R}^n}
\frac{1}{|B(x,\alpha)|}\int_{B(x,\alpha)} u^2(z)dz dx
\]
which is obtained from $L$ by averaging the original Lagrangian
over balls of radius $\alpha$. Here $B(x,\alpha)$ denotes the ball
of radius $\alpha$ about the point $x $ in $\mathbb{R}^n$ and
$| B(x,\alpha) | $ denotes its volume.

Taylor expanding the integrand about $x$ and then integrating by parts
yields cancellation of  all but the zeroth-order term, thus reproducing
exactly the original Lagrangian $L$.  This is to be expected since the
quadratic nonlinearity is rather weak, and since absolutely no information
concerning the local spatial structure of the continuum is being provided.
The latter issue is of fundamental importance and is the foundation  upon
which we shall build our theory.

\subsection{The Averaging Construction.}

To implement our construction, we will average over an ensemble
of Lagrangian fluctuation maps. We will now proceed to develop this
formalism.

\paragraph{Fuzzying the Lagrangian Flow.} Let $(M, g )$ be a $ C
^{\infty}$ compact, oriented $n$-dimensional
Riemannian manifold with $C ^{\infty}$ boundary (possibly empty).
    We consider a
two-parameter family of volume-preserving diffeomorphisms
$\xi^{\epsilon,\theta}$ of $M$ depending on  a ``radial'' component
$\epsilon \in [-\alpha/2,\alpha/2]$, $\alpha >0$,  and an ``angular''
component $\theta \in S^{n-1}_+$, where
$S^{n-1}_+$ denotes the upper hemisphere of the unit sphere
$S^{n-1}$ in ${\mathbb R}^n$.
In case $M$ has nonempty boundary, we embed $M$ into its double $\tilde{M
}$ and consider this two-parameter family defined on $\tilde{M}$; in this
case, $\xi^{\epsilon, \theta}$ need not leave $\partial M $ invariant.
This fact will be important later for certain ellipticity properties.

The parameterization is chosen such that
\begin{equation}\nonumber
\begin{array}{c}
\xi^{0,\theta}(t,x)=x, \\
\text{dist}(x,\xi^{\epsilon, \theta}(x)) < | \epsilon |
\end{array}
\end{equation}
for all $\epsilon \in [- \alpha/2,\alpha/2]$, all $t$, and $\theta \in
S^{n-1}_+$. We  define the {\bfi infinitesimal fluctuation vector} by
$$
\xi'(\theta,t,x) = \left.{\frac{d}{d \epsilon}}\right|_{\epsilon =0}
\xi^{\epsilon,\theta}(t,x),
$$
a vector field depending on the parameter $\theta$ and  time $t$.

For each time $t$, the Lagrangian flow map $\eta_t$, where $\eta_t(x) =
\eta(t,x )$, associated with a solution of the Euler equations is a
volume-preserving diffeomorphism of
$M$ which maps fluid particles $x \in M$ to
$\eta_t(x) \in M$. Motivated by the idea that a particle in a continuum
should really be regarded as a representative of a sample of particles
over a region,  we define the  $\xi ^{\epsilon, \theta} _t $-perturbed
particle placement field  by
\begin{equation}\label{1}
\eta^{\epsilon,\theta}_t(x) = (\xi^{\epsilon,\theta}_t)^{-1} \circ
\eta_t(x)
\end{equation}
for all $ \epsilon\in [- \alpha /2,\alpha/2]$ and $ \theta\in
S^{n-1}_+ $. The family of maps $\eta^{\epsilon,\theta}_t$  is called the
{\bfi fuzzy flow}.
For each $\epsilon$, $\theta$, and $t$, the map $\eta^{\epsilon,\theta}_t:
M \rightarrow M$ is a volume-preserving diffeomorphism of the fluid container.
Note that at $\epsilon=0$, $\eta^{0,\theta}_t = \eta_t$ for all
$\theta\in S^{n-1}_+$.

We take $\eta _t \in \mathcal{D}_\mu^s(M)$ and  $\xi _t ^{\epsilon, \theta}
\in \mathcal{D}_\mu^\infty( \tilde{M } )$ so that $\eta_t ^{\epsilon, \theta }
\in \mathcal{D}_\mu^s ( \tilde{M } )$. See \S\ref{review-section} for the
definition of the group $\mathcal{D}_\mu^s$.

\paragraph{Decomposition of the Spatial Velocity Field.} Our goal is to
derive the Eulerian velocity field
$u^{\epsilon,\theta}_t$ corresponding to the $\xi ^{\epsilon, \theta}
_t$-perturbed particle placement field
$\eta^{\epsilon,\theta}_t$, and define a new Lagrangian by averaging the
velocity $u^{\epsilon,\theta}_t$ over the radial parameter $\epsilon$
and the angular coordinate $\theta$.
We shall proceed with this averaging process as follows:
we begin by defining the Eulerian vector fields
associated with our three Lagrangian maps.  Let
\begin{align*}
\partial_t \eta(t,x) &= u(t,\eta(t,x)), \\
\partial_t \xi^{\epsilon,\theta}(t,x) &=
w^{\epsilon,\theta}(t,\xi^{\epsilon,\theta}(t,x)),\\
\partial_t \eta^{\epsilon,\theta}(t,x)& =
u^{\epsilon,\theta}(t,\eta^{\epsilon,\theta}(t,x)).
\end{align*}

Differentiating the Lagrangian decomposition (\ref{1}) with respect
to time $t$,
we obtain the spatial velocity decomposition
\begin{equation}\label{2}
u^{\epsilon,\theta}(t,x) =
\left\{{\xi^{\epsilon,\theta}}^*(u-w^{\epsilon,\theta})\right\}(t,x),
\end{equation}
where the notation ${\xi^{\epsilon,\theta}}^*$ denotes the pullback by
the map $\xi^{\epsilon,\theta}$.  We can also write this decomposition using
the push-forward notation via the relation
$\big( {\xi^{ \epsilon, \theta }
_t} \big) ^\ast  = {\big( \xi^{\epsilon,\theta}_t \big) ^{-1}}_*$, so that
the action on a vector field $v$ is given by
$$
\big( {\xi^{ \epsilon, \theta }
_t} \big) ^\ast v = T (\xi^{\epsilon,\theta}_t)^{-1}
\circ v \circ \xi^{\epsilon,\theta}_t,$$ where we use the symbol $T$ to denote
the tangent map (which is locally represented by the matrix of partial
derivatives).  Thus, the decomposition (\ref{2}) may be equivalently written
as
\footnote{This decomposition can also be written as
$
u^{\epsilon,\theta}(t,x) =
\text{Ad}_{(\xi_t^{\epsilon,\theta})^{-1}}(u_t-w_t^{\epsilon,\theta})$,
where $\text{Ad}$ is
the adjoint action of the volume-preserving diffeomorphism group
on diver\-gence-free vector fields.}

\begin{equation}\label{2a}
u^{\epsilon,\theta}(t,x) = T ( \xi^{\epsilon,\theta}_t)^{-1}(x)
\circ
\left( u(t, \xi^{\epsilon,\theta}_t(x)) - w^{\epsilon,\theta}(t,
\xi^{\epsilon,\theta}_t(x))\right),
\end{equation}
where, again,  $u^{\epsilon,\theta}(t,x)$ is the  Eulerian spatial velocity
field corresponding to the fuzzy flow $\eta^{\epsilon,\theta}_t$.

\paragraph{Comments on the Nature of the Decomposition.}  The Lagrangian
decomposition (\ref{1}) which ``fuzzies'' the Lagrangian flow map  yields
the decomposition (\ref{2a}) for the corresponding Eulerian variables which
is of a  {\it hybrid Lagrangian-Eulerian type}.  The Lagrangian
characteristics of this decomposition are encompassed in the presence of
the purely Lagrangian fluctuation maps $\xi^{\epsilon,\theta}_t$, and it is
indeed the presence of this Lagrangian term in (\ref{2a}) which allows
us to proceed with an asymptotic expansion which is both philosophically
and mathematically different from the ``naive'' expansion we discussed
earlier.  We should also emphasize  that without this Lagrangian aspect,
the decomposition (\ref{2a}) would reduce to the usual additive
(Reynolds) decomposition of spatial velocity fields into their mean and
fluctuating parts, which does {\it not} reflect the fuzzyness of the
Lagrangian flow map.

Our approach should also be
contrasted
with the approach taken by \cite{ChFoHoOlTiWy1998} and \cite{Holm1999}.  In
those papers, the decomposition
$$ \eta^\sigma(t,x) = \eta(t,x) + \sigma(t,\eta(t,x)),$$
is made, where $\sigma(t,x)$ is a fluctuation vector field, and
$\eta^\sigma(t,x)$ is a perturbed Lagrangian trajectory of the reference
element $x$.  This decomposition is intrinsically problematic, in that a
material vector field $\sigma(t,\eta(t,\cdot))$ is being added to a
volume-preserving diffeomorphism $\eta(t,\cdot)$.  As a consequence, the
perturbed trajectory $\eta^\sigma(t,x)$ does not come from a volume-preserving
diffeomorphism of the fluid container, that is, $\eta^\sigma(t, \cdot)$
is not a volume-preserving map.

\paragraph{The Averaged Lagrangian.} We define the {\bfi averaged
Lagrangian}
$L^\alpha$ by
\begin{align*}
L^\alpha(u) &= {\frac{1}{2}}\int_M \frac{1}{\alpha}
\int_{-\alpha/2}^{\alpha/2}  \int_{S^{n-1}_+}
\langle u^{\epsilon, \theta}(t,x), u^{\epsilon, \theta}(t,x)\rangle
d \epsilon \ \nu(\theta)\ \mu(x)  \\
&=
{\frac{1}{2}}\int_M \frac{1}{\alpha}
\int_{-\alpha/2}^{\alpha/2} \int_{S^{n-1}_+}
\left| {\xi^{\epsilon, \theta}}^*(u-w^{\epsilon, \theta})(t,x)\right|^2
d \epsilon \ \nu(\theta)\ \mu(x),
\end{align*}
where $\mu$ is the Riemannian volume form on $M$, and $\nu$ is the induced
Riemannian volume form on $S^{n-1}_+$, the upper hemisphere
of the unit sphere
$S^{n-1}$ in ${\mathbb R}^n$.

\paragraph{Comments on the Nature of the Fuzzying Operation.}
By using the upper hemisphere $S^{n-1}_+$
and integrating from $-\alpha/2$ to $\alpha/2$, we are tacitly assuming
that there is a hyperplane of symmetry in the $\theta$-parameter space.
This is not a restriction even near the boundary, as the hyperplane of
symmetry can always be chosen orthogonal to the boundary and the maps
$ \xi^{ \epsilon, \theta}_t$ can be chosen to be symmetric about
this hyperplane with respect to the radial parameter $\epsilon$.

\begin{figure}
\begin{center}
\includegraphics[scale=.7,angle=0]{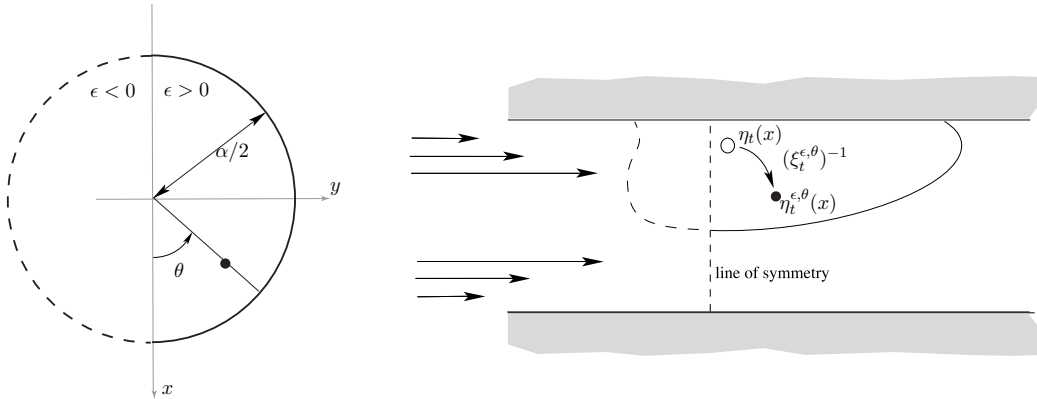}
\end{center}
\caption{\footnotesize
     An example of a perturbing map $\xi^{\epsilon, \theta}_t$; its inverse
takes the flow
point $\eta _t(x)$ to the perturbed flow point $\eta ^{\epsilon, \theta}_t
(x)$. The parameter space for $(\epsilon, \theta) $ and the symmetry
plane is shown on the left.}
\label{perturb_map.figure}
\end{figure}

The reader should keep in mind that the variables $\theta$ and $\epsilon$
parameterize possible families of maps and are not to be confused with
spatial spheres in the flow itself. We {\it are averaging over these
families of maps and not literally over spatial regions}. A representative
of the family of fluctuation maps
$\xi ^{\epsilon, \theta}$ in the two dimensional case and near a boundary is
shown in Figure \ref{perturb_map.figure}.

The internal structure behind the fuzzyness of the macroscopic
Lagrangian flow\footnote{For a general
continuum, the information about the structure of the representative sampling
would be encoded in the fluctuation maps. For example,  this might include
defects or microstructure.}
is completely encoded in the fluctuation maps
$\xi^{\epsilon,\theta}_t$.

The zeroth-order assumption that these maps
are simply the identity map leads to the {\bfi zeroth-order Lagrangian}
$L^\alpha_0$ which is exactly equal to the  Lagrangian $L$  given in
(\ref{lagrangian_euler.equation}) and thus produces  the usual Euler
equations of hydrodynamics as the continuum model.  We proceed to obtain
the first order correction to this model which accounts for the spatial
fluctuations.

\paragraph{Asymptotic Expansion.} We Taylor expand $u^{\epsilon,\theta}$ in
$\epsilon$ about
$\epsilon=0$  to obtain
\begin{equation} \label{taylor.equation}
u^{\epsilon,\theta}(t,x)=
{\xi^{\epsilon,\theta}_t}^*(u_t-w^{\epsilon,\theta}_t)(x) = u_t(x) +
\epsilon \pounds_{\xi'(\theta,t,x)} u_t(x) - \epsilon
\dot{\xi} ^\prime (\theta,t,x) + O(\epsilon^2),
\end{equation}
where the overdot means the time derivative. This  follows from the
definition of the Lie derivative, the fact that
$w^{\epsilon,\theta} = \partial_t \xi^{\epsilon,\theta}_t \circ
\xi^{\epsilon,\theta}_t$, and that
$\left. \xi^{\epsilon,\theta}_t \right| _{\epsilon = 0 } =e$.  Using the
zero-torsion condition on the Levi-Civita covariant derivative,
$\pounds_{{\xi}'} u = \nabla_{{\xi}'}u - \nabla_u {\xi}'$,  and
suppressing the dependence on $t$ and $x$, we get
$$ u^{\epsilon,\theta} = u + \epsilon( \nabla u \cdot \xi'(\theta) -
\nabla \xi'(\theta) \cdot u - \dot{\xi}'(\theta)) + O(\epsilon^2)$$
or, in index notation,
$$
{u^\epsilon }^i = u^i + \epsilon \left( u^i_{\; ;j} \xi^j(\theta) -
\xi^i_{\; ;j}(\theta)u^j - \dot{\xi}^i(\theta) \right) + O(\epsilon^2),$$
where
$$u=u^i\partial_i,\ \ \ \xi'= \xi^i \partial_i \quad \text{and}
\quad \nabla u\cdot\xi' = u^i_{\; ;j}\xi^j \partial_i.$$

In order to proceed, we make the {\bfi first-order Taylor Hypothesis} that the
infinitesimal fluctuation vector $\xi'$ is frozen, {\it as a one-form},
into the fluid so that its Lie transport vanishes;  namely,
\begin{equation}\label{3}
(\dot{\xi}'_t)^\flat + \pounds_u(\xi')^\flat = 0.
\end{equation}
We again express the Lie derivative of the $1$-form field $(\xi')^\flat $
in terms of the covariant derivative to obtain, in index notation,
$$ \dot\xi_i + u^j_{;i} \xi_j + \xi_{i;j} u^j = 0. $$
     From this hypothesis, the $O(\epsilon)$ term in the Taylor expansion
(\ref{taylor.equation}) is
\begin{align*}
u_{i;j} - \xi_{i;j} u^j + u^j_{;i} \xi_j + \xi_{i;j} u^j & =
u_{i;j} \xi^j + u^j_{;i}\xi_j \\
&= u_{i;j} \xi^j + u_{j;i}\xi^j =
2\operatorname{Def}u^\flat \cdot \xi'(\theta).
\end{align*}
It follows that
\begin{equation} \label{uet.equation}
      u ^{\epsilon, \theta} = (\xi^{\epsilon,\theta} )^\ast (u-w^\epsilon) = u
+ 2\epsilon\operatorname{Def}u \cdot \xi'(\theta) + O(\epsilon^2).
\end{equation}

Substitution of (\ref{uet.equation}) into the averaged Lagrangian
$L^\alpha$ yields
\begin{align}
& L^\alpha(u) = {\frac{1}{2}}\int_M {\frac{1}{\alpha}}
\int_{S^{n-1}_+}
\int_{-\alpha/2}^{\alpha/2}
      \left[ |u(x)|^2  + 2 \epsilon\langle u(x),
\operatorname{Def}u(x)\cdot \xi'(x,\theta)\rangle \right. \nonumber \\
& \qquad \qquad \left. + \; 4 \epsilon^2 | \operatorname{Def}u(x)\cdot
\xi'(x,\theta)|^2 + O(\epsilon^3)\right] d \epsilon \ \nu(\theta) \  \mu(x).
\label{av_Lag_approx.equation}
\end{align}

An important point about this expression is the following: There is no
contribution from the term $\langle u, O(\epsilon^2)\rangle$ to the energy
at order
$O(\epsilon^2)$. In fact, the
$O(\epsilon^2)$ term in (\ref{uet.equation}) has the form  $O
( \epsilon ^2 ) = \epsilon ^2 \left(\dot{\xi''} + R \right) $, where
\[
\dot{\xi''}:= \frac{d}{dt} \left. \frac{ d^2
\xi}{d \epsilon^2} \right| _{\epsilon=0} .
\]
However,  $\xi''$ is an independent field and must have its own dynamics
specified. We  assume that this dynamics is chosen so that $\dot{\xi}'' + R
$ is $O (\epsilon) $ and so the {\it a priori}
$O (\epsilon^2) $-term  in (\ref{uet.equation}) is in fact $O(\epsilon ^3 )$.

Integrating (\ref{av_Lag_approx.equation}) in $\epsilon$, rescaling
$\alpha \mapsto  \sqrt{\alpha/6}$, and defining
the symmetric rank-$2$ contravariant spatial fluctuation tensor (indices
up) $F$ by
$$
F(x) =
\int_{S^{n-1}_+} \xi'(x,\theta) \otimes \xi'(x,\theta) \nu(\theta),
$$
we obtain
the {\bfi first-order averaged Lagrangian}
\begin{align}
L^\alpha_1(u,F) &= {\frac{1}{2}} \int_M \int_{S^{n-1}_+}
\left[ |u|^2 + 2\alpha^2 |\operatorname{Def}u
\cdot\xi'(x,\theta) |^2 \right] \nu(\theta) \ \mu(x)\nonumber \\
&= {\frac{1}{2}} \int_M \left[ |u(x)|^2 + 2\alpha^2 \langle
F(x) \circ \operatorname{Def}u(x), \operatorname{Def}u(x) \rangle
\right] \mu(x) \nonumber\\
&= {\frac{1}{2}} \int_M \left\{ |u(x)|^2 + 2\alpha^2
\left[g(x)\otimes F(x)\right] :
\left[\operatorname{Def}u(x)\otimes \operatorname{Def}u(x)\right]
\right\} \mu(x).
\label{Lag}
\end{align}
In coordinate notation, the first-order averaged Lagrangian takes the form
$$
L^\alpha_1(u,F) = \frac{1}{2} \int_M \left\{
g_{ik}u^i u^k + 2 \alpha^2 g_{ik} F^{jl} [\operatorname{Def} u]^i_j
[\operatorname{Def} u]^k_l  \right\}
[\text{det } g]^{\frac{1}{2}} dx.
$$
The first-order averaged Lagrangian $L^\alpha_1$ is a function of
the {\bfi macroscopic Eulerian velocity field $u$} and the {\bfi
contravariant spatial fluctuation tensor $F$}.

\paragraph{The Isotropic Case.} In the case that the fluctuation tensor is
isotropic so that
$$ F(x) = g^{-1}(x),$$
the {\bfi isotropic first-order averaged Lagrangian} $L^\alpha_{1,\text{iso}}$
is given by
$$
L^\alpha_{1, {\rm iso } } (u) = {\frac{1}{2}}
\int_M \left[ |u(x)|^2 + 2\alpha^2 |\operatorname{Def}u(x)|^2 \right] \mu(x).
$$
In this special case, the Lagrangian depends only on the Eulerian velocity
field $u$ and no semi-direct product theory is required;  in fact, the
standard Euler-Poincar\'{e} theory  for reduced Lagrangian variational
principals may be invoked to obtain the {\bfi isotropic averaged Euler
equations} as
\begin{equation}\label{iso}
\begin{array}{c}
\partial_t(1-\alpha^2\Delta _r)u + \nabla_u(1-\alpha^2\Delta _r)u
-\alpha^2 (\nabla u)^t \cdot \Delta _r u = -\operatorname{grad} p,\\
\text{div } u =0, \ \ u(0)=u_0,\\
\end{array}
\end{equation}
where $\Delta _r = -(d \delta + \delta d) + 2\operatorname{Ric} $
(see \cite{Shkoller2000}).  As we stated above, the equations (\ref{iso}) are
precisely the equations of inviscid second-grade non-Newtonian fluids, and
exactly
coincide with Chorin's vortex blob algorithm when a particular choice of
smoothing kernel is used (see Oliver and \cite{Shkoller2000}). In the
case that $M$ is a manifold without boundary, the
incompressibility of the fluid allows us to replace the term
$2\alpha^2|\operatorname{Def}u|^2$ with
$\alpha^2 |\nabla u|^2$ in (\ref{Lag}),
and still obtain the identical evolution
equations as in (\ref{iso}); however, for domains with boundary it is
essential to retain the strain tensor $\operatorname{Def}u$ in the Lagrangian
so as to obtain the natural boundary conditions which ensure ellipticity
of the operator $(1-\alpha^2 \Delta _r)$.

\section{The Variational Principle and Semidirect products}

We shall next explain the sense in which the Lagrangian $L^\alpha_1(u,F)$
defined in (\ref{Lag}), a function of spatial variables, can be obtained
from a Lagrangian defined in material variables. This will be
done via an Euler-Poincar\'e procedure, which involves the group
$\mathcal{D}_{\mu, D}^s
\,\circledS\, C^\infty(T^{2,0}(M))$, the semi-direct product of the
volume-preserving diffeomorphism group $\mathcal{D}_{\mu, D}^s$ (with
Dirichlet boundary conditions) and the  smooth sections of the vector
bundle $T^{2,0}(M)$, consisting of second-rank contravariant symmetric
tensors. Before proceeding to our specific example,  we shall digress
briefly to explain the general theory.

\subsection{Lagrangian Semidirect Product Theory.}
\paragraph{The General Set Up.} Let $V$ be a
vector space and assume that the Lie group $G$ acts
linearly {\it on the right\/} on $V$ (and hence
$G$ also acts on its dual space $V^\ast$).  In the case that the
vector space $V$ consists of sections of a vector bundle $E$,
$V^*$ will denote the sections of the dual bundle $E^*$.
The semidirect
product $S = G\,\circledS\, V$ is the Cartesian product
$S  = G \times V$ whose group multiplication is given by
\begin{equation}\label{semidirectright}
(g_1, v_1) (g_2, v_2) = (g_1 g_2, v_2 + v_1 g_2),
\end{equation}
where the action of $g \in G$ on $v \in V$ is denoted
simply as $vg$.  The Lie algebra of $S$ is the
semidirect product Lie algebra,
$\mathfrak{s}    = \mathfrak{g}  \,\circledS\, V $, whose
bracket is
\begin{equation}\label{semidirectalgebraright}
[(\psi_1,v _1), (\psi_2, v_2)]
= ([\psi_1,\psi_2],\, v_1 \psi_2 - v_2 \psi_1),
\end{equation}
where we denote the induced action of $\mathfrak{g}$ on
$V$ by concatenation, as in $v_2\psi_1 $.
For $v \in V $ and $a \in V ^\ast $, define the bilinear
operator $v \diamond a  \in \mathfrak{g}^\ast$ by
\[
\left\langle v
\diamond a\,,\eta \right\rangle  = \left \langle a\eta ,
v \right \rangle.
\]

\paragraph{The Objects in Our Case.} We choose $G$ to be the topological group
$\mathcal{D}_{\mu,D}^s$. While this is not a Lie group, right multiplication
is a smooth operation, and this is the crucial feature we shall make use of.
The tangent space at the identity $T_e {\mathcal D}_{\mu,D}^s$ is equal to
$\mathfrak{X} _{{\rm div},D}^s $, the $H ^s $ vector fields on $M$ vanishing
on the boundary and with zero divergence, and plays the role of the
Lie algebra $\mathfrak{g}$.

We set $V=H^{s}(T^{2,0}(M))$, the $H^s$ sections of the vector
bundle $T^{2,0}(M)$ consisting of contravariant symmetric two
tensors (indices down). Thus, the vector space
$V ^\ast$ is $H^{s} (T^{0,2}(M))$, the $H^s$ sections of the
covariant two tensors (indices up).  The duality is with respect to
the induced Riemannian metric on $T^{2,0}(M)$.
The topological group $\mathcal{D}_{\mu,D}^s$ acts on the vector space
$H^{s} (T^{2,0}(M))$ by pull-back; hence, this action takes values in
$H^{s- 1} (T^{2,0}(M))$. Since the
group is volume preserving, the induced right action on $V^\ast$ is also by
pull-back.  We have the map $(\eta, F) \mapsto \eta^*F$.

It follows that the infinitesimal action is by the Lie derivative
which also maps $H^s$ sections into sections of class $H^{s-1}$.
Thus,
according to the above definition, the diamond operator is computed as
follows: Let $K \in H^{s-1}(T^{2,0}(M))$, $F
\in  V ^\ast = H^{s} (T^{0,2}(M))$ and let $u \in \mathfrak{g} =
\mathfrak{X} _{{\rm div},D}^s$.  We define the operator
$$\mathbf{\mathcal L}_F :\mathfrak{X} _{{\rm div},D}^s \rightarrow
H^{s-1} (T^{2,0}(M)), \ \ \
\boldsymbol{\mathcal{L}}_F (u) =
\pounds_u F.$$
Then the adjoint operator (with respect to the Riemannian metric on
$H^{s} (T^{2,0}(M))$) $\boldsymbol{\mathcal{L}}_F^*:
H^{s-1} (T^{0,2}(M)) \rightarrow \mathfrak{X} _{{\rm div},D}^s$ and is
defined by
\[
\left\langle K \diamond F, u \right\rangle =
\left\langle \pounds_u F, K \right\rangle = \left\langle u,
\boldsymbol{\mathcal{L}}_F ^\ast K \right\rangle.
\]
Thus, we have
\[
K \diamond F =  \boldsymbol{\mathcal{L}}_F ^\ast K.
\]

\paragraph{Semidirect Euler-Poincar\'e Reduction.} Assume we have a right
$G$-invariant function $
\mathcal{L}: T G
\times V ^\ast \rightarrow \mathbb{R}$. For $a_0 \in
V^\ast$, let $\mathcal{L}_{a_0} : TG
\rightarrow
\mathbb{R}$ be given by
$\mathcal{L}_{a_0}(v_g) = \mathcal{L}(v_g, a_0)$, so
$\mathcal{L}_{a_0}$ is right
invariant under the lift to $TG$ of the right action of
$G_{a_0}$ on $G$, where $G_{a_0}$ is the isotropy group
of $a_0$. Define
$L: {\mathfrak{g}} \times V^\ast \rightarrow \mathbb{R}$ by
\[
L(v_gg^{-1}, a_0g^{-1}) = \mathcal{L}(v_g, a_0).
\]

For a curve $g(t) \in G, $ let
$\xi (t) := \dot{g}(t) g(t)^{-1}$ and let $a(t) =
a_0g(t)^{-1}$, which is the unique solution of the
equation $\dot a(t) = -a(t)\xi(t)$ with initial
condition $a(0) = a_0$.

In our setting, $a_0=F_0 \in C^\infty(T^{2,0}(M))$,
\[
v_g =
u_\eta \in \{v\in H^s(M,TM) \mid v \circ \eta^{-1} \in
\mathfrak{X} _{{\rm div},D}^s \cap H ^1 _0(TM), \eta \in {\mathcal
D}_{\mu,D}^s\} ,
\]
and $L(u_\eta\circ \eta^{-1}, \eta^*F_0)={\mathcal L}(u_\eta,F_0)$
where $L$ is given by  (\ref{Lag}).

\begin{theorem} \label{ep_right}
The following are equivalent:
\begin{enumerate}
\item [{\bf i} ] Hamilton's variational principle
\begin{equation} \label{hamiltonprincipleright1}
\delta \int _{t_1} ^{t_2} \mathcal{L}_{a_0}(g(t), \dot{g} (t)) dt = 0
\end{equation}
holds, for variations $\delta g(t)$
of $ g (t) $ vanishing at the endpoints.
\item [{\bf ii}  ] $g(t)$ satisfies the Euler--Lagrange
equations for $\mathcal{L}_{a_0}$ on $G$.
\item [{\bf iii} ]  The constrained variational principle
\begin{equation} \label{variationalprincipleright1}
\delta \int _{t_1} ^{t_2}  L(\xi(t), a(t)) dt = 0
\end{equation}
holds on $\mathfrak{g} \times V ^\ast $, using variations of the form
\begin{equation} \label{variationsright1}
\delta \xi = \dot{\eta } - [\xi , \eta ], \quad
\delta a =  -a\eta ,
\end{equation}
where $\eta(t) \in \mathfrak{g}$ vanishes at the
endpoints.
\item [{\bf iv}] The Euler--Poincar\'{e} equations hold on
$\mathfrak{g} \times V^\ast$
\begin{equation} \label{eulerpoincareright1}
\frac{d}{dt} \frac{\delta L}{\delta \xi} = -
       \operatorname{ad}_{\xi}^{\ast} \frac{ \delta L}{ \delta
\xi} - \frac{\delta L}{\delta a} \diamond a.
\end{equation}
\end{enumerate}
\end{theorem}

\subsection{Computation of the Anisotropic Averaged Euler Equations}
\label{derivation_section} It is convenient to define the linear operator
${\mathcal C}:
\mathfrak{X}_{\operatorname{div}}^s
\cap H^1_0 \rightarrow H^{s-2}$, $s\ge 1$,
mapping divergence-free vector fields to vector fields, by
$$\mathcal{C}u := \operatorname{Div}\left[ C :\nabla
u^\flat\right],$$
where $\flat$ is the map from vector fields  to one-forms associated with
the metric $g$, and again the fourth-rank symmetric positive tensor $C$ is the 
symmetrization of the tensor $F\otimes g^{-1}$, given in local coordinates
by
$$ C^{ijkl}= {\frac{1}{4}}\left(
F^{lj}g^{ik}+ f^{kj}g^{il}+F^{li}g^{jk}+F^{ki}g^{jl}
\right).  $$

The functional derivatives of $L^\alpha_1$ with respect to $u$ and $F$ are
given by
\begin{equation}\nonumber
\frac{\delta L^\alpha_1}{\delta u} = (1-\alpha^2 \mathcal{C}) u
\end{equation}
and
\begin{equation}\nonumber
\frac{\delta L^\alpha_1}{\delta F} =
2 \alpha^2 \left[\operatorname{Def} u \right]^2.
\end{equation}
We can then compute that
\begin{equation}\nonumber
\frac{\delta L^\alpha_1}{\delta F} \diamond
       F =
2\alpha^2 F : \left[\nabla (\operatorname{Def} u^\flat) ^2
\right]^\sharp
   -
4\alpha^2 \operatorname{Div}\left( \left(\operatorname{Def}
u\right)^2 \cdot F \right)
\end{equation}
Letting $t=(\operatorname{Def}u)^2$, in index notation, we get
$$
\left[\frac{\delta L^\alpha_1}{\delta F} \diamond F \right]_k =
2 \alpha ^2 F^{ij} t_{ij;k} - 4 \alpha^2 \left[F^{ij} t_{kj}\right]_{;i}.
$$
Using Theorem \ref{ep_right}, we derive the following result.

\begin{theorem}
The Euler-Poincar\'e equations on Riemannian manifolds, 
associated to the Lagrangian
$L ^\alpha_1$ given by
(\ref{Lag}), are the following {\bfi anisotropic averaged Euler equations}:
\begin{align}
\partial_t (1-\alpha^2\mathcal{C})u & + \nabla_u(1-\alpha^2\mathcal{C})u
- \alpha^2
[\nabla u]^t \cdot \mathcal{C}u +
2\alpha^2 F : \left[\nabla (\operatorname{Def} u^\flat) ^2
\right]^\sharp \nonumber\\
& -
4\alpha^2 \operatorname{Div}\left( \left(\operatorname{Def}
u\right)^2 \cdot F \right)
= -\operatorname{grad} p
\label{noniso1}
\end{align}
together with the advection equation
\begin{equation}\label{noniso2}
\partial_t F + \pounds_uF =0,
\end{equation}
the incompressibility constraint
$\text{div }u=0$, initial data $u(0)=u_0$ and $F(0)=F_0$, and no-slip
conditions $u=0$ on the boundary.
\end{theorem}

\paragraph{Anisotropic Averaged Euler Equations in General Coordinates.}
In general coordinates on a manifold,
the averaged Euler equations read
\begin{align*}
& \partial_t \left(u^i - \alpha^2[C^{ijkl}u_k,_j],_l \right)
+\left( u^i - \alpha^2[C^{ijkl} u_k,_j],_l\right),_m u^m
-\alpha^2 u_m,_i[C^{mjkl}u_k,_j],_l \\
& \qquad  + 2\alpha^2 F^{kj} [D_{km}g^{mn}
D_{nj}],_i
      -4\alpha^2[F^{kj} D_{im}g^{mn} D_{nj}],_k   \\
& \qquad \quad = -p,_i
\end{align*}
where, as earlier, $D _{ij} = \frac{1}{2}\left(u_{i;j} + u_{j;i}\right)$
is the rate of deformation tensor and indices are raised and lowered using
the metric tensor (which of course need not be diagonal in general
coordinates),  and
$ C^{ijkl}= {\frac{1}{4}}\left(
F^{lj}g^{ik}+ f^{kj}g^{il}+F^{li}g^{jk}+F^{ki}g^{jl}
\right)$.
 In Euclidean space, one need only set the components of
the metric tensor $g_{ij}$ to the Kronecker delta $\delta_{ij}$.

\paragraph{Comments on the Form of the Equations.} In 2D, identifying $F$
with the vector
$(F^{11}, F^{12}, F^{22})$, equation (\ref{noniso2}) takes the form
           \begin{equation}
           \frac{D}{dt}
           \left[\begin{array}{cccc}
            F^{11}\\
            F^{12} \\
            F^{22}
           \end{array} \right]
           =
           \left[\begin{array}{cccc}
           2u^1_{,1} & 2u^1_{,2} &  0         \\
           u^2_{,1}  & 0         &  u^1_{,2}  \\
           0         & 2u^2_{,1} & -2u^1_{,1}  \\
           \end{array} \right]
           \left[\begin{array}{cccc}
           F^{11} \\
           F^{12} \\
           F^{22} \\
           \end{array} \right],
           \label{2D}
           \end{equation}
where $D/dt$ denotes $\partial_t + (u \cdot \nabla)$.
Notice that the matrix on the right-hand-side of (\ref{2D}) is
traceless; a similar form holds in 3D as well.  This is not surprising,
since by virtue of the incompressibility of the Lagrangian flow and
the fact that $F_t = \eta_t^* F_0$, we have
that
\[
\det(F_t) = \det(F_0),
\]
for all $t$ for which the solution exists.  As consequence, the operator
$(1-\alpha^2 {\mathcal C})$ remains uniformly elliptic, if $F_0$ is strictly
positive.

\paragraph{The Circulation Theorem.}
Let $\gamma : S^1 \rightarrow M $ be a loop and let $\gamma _t =
\eta _t \circ \gamma$ denote the evolution of the loop moving with the
fluid.

\begin{theorem} For a solution of the anisotropic averaged Euler equations,
we have
\[
\frac{d}{dt} \int _{\gamma _t}  (1 - \alpha ^2 \mathcal{C}^\flat ) u
^\flat = 2 \alpha ^2 \int _{\gamma _t}  \left[
2\operatorname{Div}\left( \left(\operatorname{Def}
u\right)^2 \cdot F \right)-
   F : \left[\nabla (\operatorname{Def} u^\flat) ^2
\right]^\sharp
    \right] ^\flat.
\]
\end{theorem}

This follows directly from the Lie derivative form of the equations given
in (\ref{Lie_average.equation}). We note that if one were to make use of the
general Kelvin-Noether theorem given in \cite{HoMaRa1998b}, one would
arrive at the same result.

\section{Analytic Properties}

In this section we prove well-posedness and other properties of the
solutions by showing that these equations are given by a smooth vector
field in material representation in the appropriate Sobolev topologies.
This is in line with what is known about the Euler equations, as described
in the introduction. We also discuss the corrector for the equations.

\subsection{Well-posedness of Classical Solutions}

We shall prove existence, uniqueness, and {\it smooth} dependence on initial
data on finite time intervals
for solutions of the anisotropic averaged Euler equations.  For simplicity,
we shall restrict the fluid domain $M$ to be a compact subset of Euclidean
space with smooth boundary,
although our methods can be applied to Riemannian manifolds.

We begin by collecting some preliminary results.  Set
${\mathcal V}^s = H^1_0 \cap H^s$ and
${\mathcal V}^s_\mu =
H^1_0 \cap \mathfrak{X}_{\operatorname{div}}^s$. Also, let
${\mathcal D}_{D}^s$
denote the $H^s$ class diffeomorphisms which fix the boundary, and
again let
${\mathcal D}_{\mu,D}^s$ denote the diffeomorphisms in ${\mathcal D}_{D}^s$
which preserve the volume $\mu$.

\begin{lemma} \label{lemma1}
For $u\in {\mathcal V}^s_\mu$, $s>1$,
$$
\partial_t \mathcal{C}u = \mathcal{C}(\partial_t u)
+\operatorname{Div} \left[-\nabla _u F \cdot \nabla u + \nabla u\cdot
\nabla u \cdot
F + \nabla u \cdot F \cdot \nabla u^t\right],
$$
and
\begin{align*}
\nabla_u \mathcal{C}u &= \mathcal{C}(\nabla_uu) +
\nabla u \cdot \nabla_u\operatorname{Div} F + \nabla \nabla u : \nabla_uF
-2 \nabla\nabla u :(\nabla u \cdot F)\\
&\qquad - \nabla u \cdot(\nabla u \cdot
\operatorname{Div} F) - \nabla u \cdot (\nabla \nabla u : F).
\end{align*}
\end{lemma}
\begin{proof}
The proof is a simple computation which we leave to the interested reader,
c.f. Lemma 3 in \cite{Shkoller2000}.
\end{proof}

Set $\mathcal{L} = \operatorname{Def} ^\ast \left[ \left( g \otimes F
\right) : \operatorname{Def} \right] $. Then
$\mathcal{L}$ is a positive unbounded self-adjoint operator on $L^2$ with
domain ${\mathcal V}^2_\mu$.  Define the
inner-product $(\cdot, \cdot)$ on ${\mathcal V}^2_\mu$ by
$$(u,v) = \langle (1-\alpha^2\mathcal{L}) u, v \rangle_{L^2}.$$
For $s>n/2+1$, $(\cdot, \cdot)$ defines an inner-product on
$T_e\mathcal{D}_{\mu,D}^s$, the tangent space at the identity of the
subgroup $\mathcal{D}_{\mu,D}^s$ consisting of those elements of
$\mathcal{D}_\mu^s$ which restrict to the identity on ${\partial M}$.
Right-translating $(\cdot, \cdot)$ to the entire group $\mathcal{D}_{\mu,D}^s$
defines a $C^\infty$ weak Riemannian metric  by Proposition 3 of
\cite{Shkoller2000}.

\begin{proposition}\label{P}
For $r\ge 1$ we have the following well-defined decomposition
\begin{equation}\label{p1}
{\mathcal V}^r
= {\mathcal V}^r_\mu \oplus (1-\mathcal{L})^{-1}\operatorname{grad}
H^{r-1}(M).
\end{equation}
Thus, if $F\in {\mathcal V}^r$, then there exists $(v,p) \in
{\mathcal V}^r_\mu \times H^{r-1}(M)/{\mathbb R}$ such that
$$
F= v + (1-\mathcal{L})^{-1} \operatorname{grad}p
$$
and the pair $(v,p)$ are solutions of the Stokes problem
\begin{equation}\label{p2}
\begin{array}{c}
(1-\mathcal{L}) v + \operatorname{grad}p = (1-\mathcal{L})F, \\
\qquad \operatorname{div} v =0,\\
\qquad \qquad \qquad v = \text{\rm on } {\partial M}.\\
\end{array}
\end{equation}
The summands in (\ref{p1}) are $( \cdot, \cdot )$-orthogonal.
Now, define the Stokes projector
\begin{equation}\label{p3}
\begin{array}{c}
{\mathcal P}_e: {\mathcal V}^r \rightarrow {\mathcal V}^r_\mu,\\
{\mathcal P}_e(F) = F - (1-\mathcal{L})^{-1}\operatorname{grad}p.
\end{array}
\end{equation}
Then, for $s>(n/2)+1$, $\overline {\mathcal P}: T{\mathcal D}_{D}^s \rightarrow
T{\mathcal D}_{\mu,D}^s$, given on each fiber by
\begin{equation}\nonumber
\begin{array}{c}
\overline{\mathcal P}_\eta: T_\eta{\mathcal D}_{D}^s \rightarrow
T_\eta{\mathcal D}_{\mu,D}^s,\\
\overline{{\mathcal P}}_\eta(X_\eta) = \left[ {\mathcal P}_e(X_\eta\circ
\eta^{-1})\right] \circ\eta,
\end{array}
\end{equation}
is a $C^\infty$ bundle map covering the identity.
\end{proposition}
\begin{proof}
The proof is identical to the proof of Proposition 2 in \cite{Shkoller2000}.
\end{proof}
\begin{theorem}\label{Thm1}
Set $s>(n/2)+2$,  and
let $\langle\langle \cdot, \cdot \rangle\rangle$ denote the right invariant
metric
on $\mathcal{D}_{\mu,D}^s$ given at the identity by $(\cdot, \cdot)$.
For $u_0 \in T_e \mathcal{D}_{\mu,D}^s$ and $F_0 \in C^\infty(
T^{2,0})$,  there exists an
interval $I=(-T,T)$, depending on $|u_0|_s$, and
a unique geodesic $\dot \eta$ of $\langle\langle \cdot, \cdot \rangle\rangle$
with initial data
$\eta(0) = e$ and $\dot \eta(0)=u_0$ such that
$$\dot \eta \text{ is in } C^\infty( I, T\mathcal{D}^s_{\mu,D})$$
and has $C^\infty$  dependence on the initial velocity $u_0$.

The geodesic $\eta$ is the Lagrangian flow of the time-dependent vector
field $u(t,x)$ given by
$$ \partial_t\eta(t,x) = u(t,\eta(t,x)),$$
and, with $F(t,x) = \left( \eta _t \right)_\ast F _0 (x)$,
$$(u,F) \in  C^0( I, {\mathcal V}^s_\mu) \cap C^1(I,{\mathcal V}^{s-1}_\mu)
\times C^0(I, H^{ s - 1 }(T^{2,0}) ).  $$
uniquely solves the anisotropic averaged Euler equations with Dirichlet
boundary conditions $u=0$, and depends continuously
on $(u_0,F_0)$.
\end{theorem}
\begin{proof}
The key to the proof rests in the fact that the pair $(u,F)$ solves the
anisotropic averaged Euler equations if and only if $\eta$ is a solution
of
\begin{equation}\label{lag}
\ddot \eta + \mathcal{U}^\alpha(\dot \eta \circ \eta^{-1}) \circ \eta
= \left[(1-\alpha^2 \mathcal{L})^{-1} \operatorname{grad} p\right] \circ
\eta,
\end{equation}
where
\begin{align*}
\mathcal{U}^\alpha(u) =& \alpha^2(1-\alpha^2\mathcal{L})^{-1}
\left\{ \operatorname{Div}
        \left[ \nabla_uF \cdot \nabla u -\nabla u\cdot\nabla u \cdot F
- \nabla u \cdot F \cdot \nabla u^t\right]\right. \\
&- \nabla u \cdot(\nabla_u\operatorname{Div}F) - \nabla \nabla u : \nabla_u
F  + 2 \nabla \nabla u:(\nabla u \cdot F) \\
&+ \nabla u \cdot(\nabla u \cdot \operatorname{Div}F)
+ \nabla u \cdot(\nabla\nabla
u : F) - \nabla u^t \cdot {\mathcal C}u\\
&\left.+ 2F:\nabla(\operatorname{Def} u^2)
         -4\operatorname{Div}[F\cdot\operatorname{Def} u^2]\right\}.
\end{align*}
This expression is obtained using Lemma \ref{lemma1}.  Now it is clear
that $\mathcal{U}^\alpha$ maps $H^s$ vector fields into $H^s$ vector fields
since $H^{s-2}$ forms a multiplicative algebra, and since the
fluctuation tensor
$F$ at $t=0$ is given by $F_0$ which is $C^\infty$.  In particular, in the
Lagrangian frame, $F$ is frozen, so the elliptic operator
$[(1-\alpha^2 \Delta )(\dot\eta\circ \eta^{-1})]\circ \eta$ has $C^\infty$
coefficients.

Thus, the proof of Theorem 2 in \cite{Shkoller2000} gives a unique curve
$\eta \in C^\infty(I, \mathcal{D}_{\mu,D}^s)$ solving (\ref{lag}).
That $u$ is only $C^0$ in time follows from the fact that the
map $\eta \mapsto \eta^{-1}:\mathcal{D}_{\mu,D}^s \rightarrow
\mathcal{D}_{\mu,D}^s$ is only $C^0$.  That $F$ is in $C^0(I,
H^{ s - 1 }(T^{2,0}))$ follows from the regularity of $u$.
\end{proof}

\subsection{A Corrector for the Macroscopic
Velocity}\label{corrector_subsection}

The solution to the anisotropic averaged Euler equations (\ref{noniso1})
and (\ref{noniso2}) yields the pair $(u,F)$. The
macroscopic spatial velocity field $u$ is only the zeroth-order term in
the expansion (in $\epsilon$) for the velocity field $u^{\epsilon,\theta}$.
We have computed, in equation (\ref{uet.equation}), the expansion of
$u^{\epsilon,\theta}$ to order $O(\epsilon^2)$ as
$$u^{\epsilon,\theta}(t,x)=
u(t,x)+2\epsilon\operatorname{Def}(u)(t,x)\cdot \xi'(t,x,
\theta) + O(\epsilon^2).$$
Since $|\epsilon|$ is bounded by $\alpha/2$, we have that
$$u^{\alpha,\theta}(t,x)=
u(t,x)+\alpha\operatorname{Def}(u)(t,x)\cdot \xi'(t,x,
\theta) + O(\alpha^2),$$
so that we may add the $O(\alpha)$ term to the expansion by solving
for the infinitesimal fluctuation vector $\xi'(t,x,\theta)$.
This, however, only requires the solution of the simple linear advection
problem (\ref{3}) given by
$$
\dot{\xi}'(t,x,\theta)^\flat + \pounds_{u(t,x)}\xi'(t,x,\theta)^\flat = 0.
$$

Computationally, this means that we may solve for the macroscopic velocity
field $u$ at spatial scales larger than $\alpha$ and correct for the
unresolved small scales to $O(\alpha^2)$.

\subsection{Limits of Zero Viscosity}  \cite{Pe1985} showed that
by perturbing the Euler solution's  Lagrangian particle trajectory
$\eta_t(x)$ by Brownian motions and averaging over such perturbations,
the Navier-Stokes equations are obtained.  In other words, letting Euler
trajectories take random-walks produces the viscosity term $\nu \Delta u$,
where $\eta_t$ is the flow map for the velocity field $u$.  In the
setting of the averaged Euler equations, the Lagrangian trajectory $\eta_t(x)$
of a particle $x$ corresponds to the flow of the velocity  $u(t,x)$ solving
the anisotropic averaged Euler equations.   Thus, Peskin's argument can
be carried over in this setting to obtain the same viscous term $\nu \Delta u$.

We are hence motivated to define the {\bfi anisotropic averaged Navier-Stokes
equations} by
\begin{align}
\partial_t (1-\alpha^2\mathcal{C})u^\nu & + \nabla_{u^\nu}
(1-\alpha^2\mathcal{C})u^\nu- \alpha^2
[\nabla u^\nu]^t \cdot \mathcal{C}u^\nu +
2\alpha^2 F : \left[\nabla (\operatorname{Def} {u^\nu}^\flat) ^2
\right]^\sharp \nonumber\\
& -
4\alpha^2 \operatorname{Div}\left( \left(\operatorname{Def}
u^\nu\right)^2 \cdot F \right)
= -\operatorname{grad} p + \nu \Delta u^\nu, \ \ \nu >0
\label{avg_NS}
\end{align}
together with the advection for the fluctuation tensor $F$ given by
(\ref{noniso2}), the incompressibility constraint $\operatorname{div}u=0$,
initial data $u(0)=u_0$ and $F(0)=F_0$, and the no-slip conditions $u=0$ on
the boundary.

Let $\eta^\nu_t$ denote the Lagrangian flow of the solution $u^\nu$ of
the anisotropic averaged Navier-Stokes equations (\ref{avg_NS}), and
let $\dot\eta^\nu$ denote the partial time derivative of the flow, i.e., the
material velocity field.

\begin{theorem}
For $s>(n/2)+2$ and $(u_0, F_0) \in {\mathcal V}^s_\mu\times
C^\infty(T^{2,0})$,
there exists a $T>0$, depending only on $\|u_0\|_{H^s}$ on not on the
viscosity $\nu$, such that for each $\nu >0$
$$\dot\eta \text{\rm \ \  is in } C^\infty ([0,T), T{\mathcal D}_{\mu,D}^s),$$
and has $C^\infty$ dependence on the initial velocity field $u_0$.
Furthermore, $u_t^\nu = \dot\eta_t^\nu \circ {\eta^\nu}_t^{-1}$ is in
$C^0( [0,T), {\mathcal V}^s_\mu) \cap C^r([0,T),{\mathcal V}^{s-r}_\mu)$
and depends continuously on $u_0$.
\end{theorem}

The proof follows the proof of Theorem 2 in \cite{Shkoller2000}; we refer the
interested reader there for the details.  As a consequence of the time
interval $[0,T)$ of solutions $u^\nu$ being independent of $\nu$, we
immediately obtain the following.

\begin{corollary}
For $s>(n/2)+2$, solutions $u^\nu$ of (\ref{avg_NS}) converge regularly to
the inviscid solutions $u$ of (\ref{noniso1}) as $\nu \rightarrow 0$.
Furthermore, letting
$u^\nu = \partial_t \eta_\nu \circ \eta_\nu^{-1}$,  the viscous Lagrangian flow
$\eta_\nu$ converges regularly in the $H^s$ topology to the inviscid
Lagrangian flow $\eta = \eta^0$.
\end{corollary}
This result states that we can generate smooth-in-time classical solutions to
the anisotropic
averaged Euler equations by obtaining a sequence of viscous solutions and
allowing $\nu$ to go to zero, and what is surprizing, this holds even in
the presense of boundaries.  Results of this type were conjectures in
\cite{MaEbFi1972} and \cite{BaCh1998a} (see also \cite{BaCh1998b}).

In the case of the isotropic averaged Euler
equations, \cite{FoHoTi2000} have added the dissipative term
$\nu \Delta (1-\alpha^2 \Delta) u$ instead of using $\nu \Delta u$, and this
is enough to give global in time classical solutions in dimension three.
It is, however, the term $\nu \Delta u$ that arises from either the approach
of \cite{Pe1985}  noted above, or from the constitutive theory approach of
\cite{RiEr1955}.

\subsection*{Future Directions} There are several interesting directions

\begin{enumerate}
\item Of course numerical implementation for specific
flows will be of great interest.
\item Modeling the mean velocity profile for turbulent flows in channels
and pipes.
\item Specific flows and special solutions.
\item Links with elliptical vortex blob methods of Zabusky and coworkers
(see, e.g., \cite{MeZaMc1988}) would be of
interest to establish; it is reasonable to expect that solutions of this
sort would provide the anisotropic analog of the vortex blob solutions of
\cite{OlSh2000}.
\item Further investigation of the vorticity formulation and its relation
with the coadjoint orbit structure in the semidirect product for the
Hamiltonian version of this theory.
\end{enumerate}

\paragraph{Acknowledgments.}
We thank G.I. Barenblatt for useful comments regarding his work on
averaging for damage propagation and the connection to our work,
Alexandre Chorin for useful discussions on links between our work and
optimal prediction methods, and Albert Fannjiang on extensions to
stochastic perturbation methods.
We also thank Tudor Ratiu for his helpful comments. In particular,
the important idea of advecting the fluctuations as a one-form rather than
as a vector field (see equation (\ref{3})) was done in collaboration with
him. We  thank Darryl Holm for keeping us informed about his work on
fluctuation effects.

JEM and SS  were partially supported by the NSF-KDI grant
ATM-98-73133. SS was partially supported by the Alfred P. Sloan Foundation
Research  Fellowship.

{\small
\bibliographystyle{new}    

\begin{thebibliography}{300}
\expandafter\ifx\csname natexlab\endcsname\relax\def\natexlab#1{#1}\fi

\bibitem[Arnold(1966)]{Arnold1966} Arnold, V.~I. [1966], Sur la
g\'{e}om\'{e}trie differentielle des groupes de {L}ie de dimension infinie
et ses applications \`{a} l'hydrodynamique des fluids parfaits,  {\it Ann.
Inst. Fourier, Grenoble\/},  {\bf 16}, 319--361.

\bibitem[Arnold and Khesin(1998)]{ArKh1998} Arnold, V.~I. and B.~Khesin
[1998], {\it Topological Methods in Hydrodynamics}, Appl. Math. Sciences,
{\bf 125},  Springer-Verlag.

\bibitem[Barenblatt and Chorin(1998a)]{BaCh1998a} Barenblatt, G.~I and
Chorin, A.~J. [1998], Scaling laws and vanishing viscosity limits in
turbulence theory,  {\it Proc. Sympos. Appl. Math.\/},   {\bf 54}, 1--25.

\bibitem[Barenblatt and Chorin(1998b)]{BaCh1998b} Barenblatt, G.~I and
Chorin, A.~J. [1998], New perspectives in turbulence: scaling laws,
asymptotics, and intermittency,  {\it SIAM Rev.\/},   {\bf 40}, 265--291,
SIAM, Philadelphia, PA.

\bibitem[Barenblatt and Prostokishin(1993)]{BaPr1993} Barenblatt, G. I. and
V.M. Prostokishin [1993], A mathematical model of damage accumulation
taking into account microstructural effects,  {\it European J. Appl.
Math.\/}, 225--240.

\bibitem[Chen et~al. (1998)]{ChFoHoOlTiWy1998}
Chen, S.~Y., C.~Foias, D.~D.~Holm, E.~J.~Olson, E.~S.~Titi and S.~Wynne
[1998], A connection between the Camassa-Holm equations and turbulent flows
in channels and pipes,  {\it Phys. Fluids\/},   {\bf 11}, 2343--2353.

\bibitem[Chorin et~al. (1999)Chorin, Kast and Kupferman]{ChKaKu1999}
Chorin, A. J., Kast, A. P. and Kupferman, R. [1999], Unresolved computation
and optimal predictions,  {\it Comm. Pure Appl. Math.\/},   {\bf 52},
1231--1254.

\bibitem[Ebin and Marsden(1970)]{EbMa1970} Ebin, D.~G. and J.~E.~Marsden
[1970], Groups of diffeomorphisms and the motion of an incompressible
fluid,  {\it Ann. of Math.\/},   {\bf 92}, 102--163.

\bibitem[Holm(1999)]{Holm1999} Holm, D.~D. [1999], Fluctuation effects on
3D Lagrangian mean and Eulerian mean fluid motion,  {\it Physica D\/},
{\bf 133}, 215--269.

\bibitem[Holm et~al. (1998a)Holm, Marsden and Ratiu]{HoMaRa1998a} Holm,
D.~D., J.~E.~Marsden and T.~S.~Ratiu [1998a], {E}uler--{P}oincar\'e models
of ideal fluids with nonlinear dispersion,  {\it Phys. Rev. Lett.\/},
{\bf 349}, 4173--4177.

\bibitem[Holm et~al. (1998b)Holm, Marsden and Ratiu]{HoMaRa1998b} Holm,
D.~D., J.~E.~Marsden and T.~S.~Ratiu [1998b], The {E}uler--{P}oincar\'{e}
equations and semidirect products with applications to continuum theories,
{\it Adv. in Math.\/},   {\bf 137}, 1--8.

\bibitem[Foias, Holm, and Titi (2000)]{FoHoTi2000} Foias, C., D.~D.~ Holm,
and E.~S.~Titi [2000], preprint.

\bibitem[Melander et~al. (1988)Melander, Zabusky and
McWilliams]{MeZaMc1988} Melander, M. V., N.J. Zabusky and J.C. McWilliams
[1988], Symmetric vortex merger in two dimensions: causes and conditions,
{\it J. Fluid Mech.\/},   {\bf 195}, 303--340.

\bibitem[Marchioro and Pulvirenti(1994)]{MaPu1994} Marchioro, C. and M.
Pulvirenti [1994], Mathematical Theory of Incompressible Nonviscous
Fluids,  Springer.

\bibitem[Marsden et~al. (1972)Marsden, Ebin and Fischer]{MaEbFi1972}
Marsden, J.~E., D.~G.~Ebin and A.~Fischer [1972], Diffeomorphism groups,
hydrodynamics and relativity, in {\it Proceedings 13th Biennial Seminar on
Canadian Mathematics Congress\/}, 135--279.

\bibitem[Marsden and Ratiu(1999)]{MaRa1999} Marsden, J.~E. and T.~S.~Ratiu
[1999], Introduction to Mechanics and Symmetry, Texts in Applied
Mathematics,   {\bf 17},  Springer-Verlag; \newblock Second Edition, 1999.

\bibitem[Marsden et~al. (2000)Marsden, Ratiu and Shkoller]{MaRaSh2000}
Marsden, J.~E., T.~Ratiu and S.~\cite{Shkoller2000}, The geometry and analysis
of the averaged Euler equations and a new diffeomorphism group,  {\it Geom.
Funct. Anal.\/}; \newblock (to appear).

\bibitem[Marsden and Weinstein(1983)]{MaWe1983} Marsden, J.~E. and
A.~Weinstein [1983], Coadjoint orbits, vortices and {C}lebsch variables for
incompressible fluids,  {\it Physica D\/},   {\bf 7}, 305--323.

\bibitem[Oliver and Shkoller(2000)]{OlSh2000} Oliver, M. and S.
\cite{Shkoller2000}, The vortex blob method as a second-grade
non-Newtonian fluid;
\newblock E-print,
\textcolor{blue}{\url{http://xyz.lanl.gov/abs/math.AP/9910088/}}.

\bibitem[Peskin(1985)]{Pe1985} Peskin, C.S. [1985]
A randon-walk interpretation of the incompressible Navier-Stokes equations,
{\it Comm. Pure Appl. Math.\/} {\bf 38}, 845--852.

\bibitem[Rivlin and Erickson(1955)]{RiEr1955} Rivlin, R.S. and J.L.
Erickson [1955], Stress-deformation relations for isotropic materials,
{\it J. Rat. Mech. Anal.\/},   {\bf 4}, 323--425.

\bibitem[Shkoller(1998)]{Shkoller1998} Shkoller, S. [1998], Geometry and
curvature of diffeomorphism groups with $H^1$ metric and mean
hydrodynamics,  {\it J. Funct. Anal.\/},   {\bf 160}, 337--365.

\bibitem[Shkoller(2000)]{Shkoller2000} Shkoller, S. [2000], On averaged
incompressible Lagrangian  hydrodynamics; \newblock E-print,
\textcolor{blue}{\url{http://xyz.lanl.gov/abs/math.AP/9908109/}}.

\end{thebibliography}

 }

\end{document}